\documentclass[letterpaper,9 pt, technote]{IEEEtran}
\IEEEoverridecommandlockouts
\usepackage[caption=false]{subfig}
\usepackage[dvips]{graphicx}
\graphicspath{{figures/}}
\usepackage{amsmath, amssymb}
\usepackage{tikz} \usetikzlibrary{calc}
\usepackage{psfrag}
\usepackage{pgfplots}
\usepackage{algorithm}     % <-- should be replaced with 'algorithms'
\usepackage{algpseudocode} % <-- should be replaced with 'algorithms'

\DeclareMathOperator*{\minimize}{minimize}

\DeclareMathOperator*{\trace}{tr}
\DeclareMathOperator*{\tr}{tr}
\DeclareMathOperator*{\supremum}{sup}

\DeclareMathOperator*{\eigen}{eig}

\DeclareMathOperator*{\subject}{subject\ to}

\DeclareMathOperator*{\diag}{diag}

\DeclareMathOperator*{\adj}{adj}

\DeclareMathOperator*{\IQC}{IQC}

\DeclareMathOperator*{\degree}{deg}
\DeclareMathOperator*{\real}{Re}

\newcounter{thm}
\newcounter{cor}

\newcounter{remcount}

\newtheorem{theorem}[thm]{Theorem}
\newtheorem{corollary}[cor]{Corollary}
\newtheorem{rem}[remcount]{Remark}

\ifCLASSINFOpdf
\else
\fi

\begin{document}
\title{Robust Stability Analysis of Sparsely Interconnected Uncertain Systems$^{\ast}$}

\author{Martin S. Andersen$^{1}$, Sina  Khoshfetrat Pakazad$^{2}$, Anders Hansson$^{2}$,  and Anders Rantzer$^{3}$% <-this % stops a space
\thanks{$^*$This work has been supported by the Swedish Department of Education within the ELLIIT project.}% <-this % stops a space
\thanks{$^{1}$Martin S. Andersen is with the Department of Applied Mathematics and Computer Science, Technical University of Denmark. Email:
        {\tt\footnotesize mskan@dtu.dk}. The work was carried out
      while Martin S. Andersen was a postdoc at Link\"oping University.}
\thanks{$^{2}$ Sina Khoshfetrat Pakazad and Anders Hansson are with the
    Division of Automatic Control, Department of Electrical
    Engineering, Link\"oping University, Sweden. Email:
        {\tt\footnotesize \{sina.kh.pa, hansson\}@isy.liu.se}. The work was carried out
      while Anders Hansson was a visiting professor at the University of
      California, Los Angeles.}%
\thanks{$^{3}$Anders Rantzer is with the Department of Automatic Control,
    Lund University, Sweden. Email: {\tt\footnotesize anders.rantzer@control.lth.se}}}% <-this % stops a space

\maketitle

\begin{abstract}
In this paper, we consider robust stability analysis of large-scale sparsely interconnected uncertain systems. By modeling the interconnections among the subsystems with integral quadratic constraints, we show that robust stability analysis of such systems can be performed by solving a set of sparse linear matrix inequalities. We also show that a sparse formulation of the analysis problem is equivalent to the classical formulation of the robustness analysis problem and hence does not introduce any additional conservativeness. The sparse formulation of the analysis problem allows us to apply methods that rely on efficient sparse factorization techniques, and our numerical results illustrate the effectiveness of this approach compared to methods that are based on the standard formulation of the analysis problem.
\end{abstract}


\IEEEpeerreviewmaketitle

\section{Introduction}\label{sec:Intro}

\IEEEPARstart{R}{obust} stability of uncertain systems can be analyzed using different approaches, e.g., $\mu$-analysis and IQC analysis \cite{robustandoptimal,multivariable,ulfiqc,rantzer}. In general, the computational burden of performing robust stability analysis grows rapidly with the state dimension of the uncertain system. For instance, performing robust stability analysis using integral quadratic constraints (IQCs) involves finding a solution to a frequency-dependent semi-infinite linear matrix inequality (LMI). This frequency dependent semi-infinite LMI can be reformulated using the Kalman-Yakubovich-Popov (KYP) lemma as a finite-dimensional frequency independent LMI which is generally dense \cite{rantzer,ran:96}. The size and number of variables of this LMI grows with the number of states and dimension of the uncertainty in the system, and hence IQC analysis of uncertain systems with high state and uncertainty dimensions can be very computationally demanding. This is the case even if the underlying structure in the resulting LMI is exploited, \cite{ChapterBook1,Anders,Parrilo,Kao2,and+dah+van10,fuj+kim+koj+oka+yam09,WHJ:09,Wallin,Kao1}. An alternative approach to solving the semi-infinite LMI is to perform frequency gridding where the feasibility of the frequency dependent semi-infinite LMI is verified for a finite number of frequencies. Like the KYP lemma-based formulation, this method leads to a set of LMIs that are generally dense and costly to solve for problems with a high state dimension. The focus of this paper is the analysis of large-scale sparsely interconnected uncertain systems, which inherently have high number of states. As a result, the mentioned methods are prohibitively costly for analyzing large-scale interconnected uncertain systems even when the interconnections are sparse. In order to address the robustness analysis of such systems, we will show how the sparsity in the interconnection can be utilized to reduce the computational burden of solving the analysis problem. Sparsely interconnected uncertain systems are common in e.g.\ large power networks that include many interconnected components, each of which is connected only to a small number neighboring components.

In \cite{Langbort} and \cite{P5}, the authors consider robust stability analysis of interconnected systems where the uncertainty lies in the interconnections among the subsystems. In these papers, the authors consider the use of IQCs to describe the uncertain interconnections, and they provide coupled LMIs to address the stability analysis and control design for such systems. However, they do not investigate the computational complexity of analyzing large-scale systems. In \cite{UlfLetter}, a method is proposed for robust stability analysis of interconnected uncertain systems using IQCs. The authors show that when the interconnection matrix is given as $\Gamma = \bar \Gamma \otimes I$ and the adjacency matrix of the network $\bar \Gamma$ is normal, the analysis problem can be decomposed into smaller problems that can be solved more easily. In \cite{Kao} and \cite{Ulf}, the authors describe a method for robust stability analysis of interconnected systems with uncertain interconnections. The analysis approach considered in these papers is based on separation of subsystem frequency responses and eigenvalues of the interconnection and adjacency matrices. Although the computational complexity of this approach scales linearly with the number of subsystems in the network, the proposed method can only be used for analyzing interconnected systems with specific interconnection descriptions. In \cite{Kao}, the uncertain interconnection matrix is assumed to be normal, and its spectrum is characterized using quadratic inequalities. In \cite{Ulf},  the interconnection matrix is defined as $\Gamma = \bar \Gamma \otimes I$, where the adjacency matrix of the network $\bar \Gamma$ is assumed to be normal with its spectrum expressed using general polyhedral constraints.  The authors in \cite{Vinnicombe} provide a scalable method for analyzing robust stability of interconnections of SISO subsystems over arbitrary graphs. The proposed analysis approach in \cite{Vinnicombe} is based on Nyquist-like conditions which depend on the dynamics of the individual subsystems and their neighbors. Finally, the authors in \cite{kim:12} propose a decomposable approach to robust stability analysis of interconnected uncertain systems. When the system transfer matrices for all subsystems are the same, this approach makes it possible to solve the analysis problem by studying the structured singular values of the individual subsystems and the eigenvalues of the interconnection matrix.   

In related work \cite{and+han:12}, we considered robust stability analysis of interconnected uncertain systems where the only assumption was the sparsity of the interconnection matrix, and we laid out the basic ideas of how to exploit this kind of sparsity in the analysis problem. In that paper, we put forth an IQC-based analysis method for analyzing such systems, and we showed that by characterizing the interconnections among subsystems using IQCs, the sparsity in the interconnections can be reflected in the resulting semi-infinite LMI. No numerical results were presented in \cite{and+han:12}. In this paper, we present an extended version of \cite{and+han:12} where we show how the sparsity in the problem allows us to use highly effective sparse SDP solvers, \cite{BeY:05,and+dah+van10,fukuda_exploitingsparsity}, to solve the analysis problem in a centralized manner. Note that none of the methods mentioned in the previous paragraph consider as general a setup as the one presented in this paper.
\vspace*{-5pt}
\subsection{Outline}
This paper is organized as follows. Section \ref{sec:IQC} briefly reviews relevant theory from the IQC analysis framework. In Section~\ref{sec:Inter}, we provide a description of the interconnected uncertain system that is later used in lumped and sparse formulations of the analysis problem. We discuss these formulations in sections \ref{sec:Lumped} and~\ref{sec:Sparse}, and we report the numerical results in Section \ref{sec:Results}. Section~\ref{sec:Conclusions} presents conclusions and some remarks regarding future research directions. 
\vspace*{-5pt}
\subsection{Notation}
We denote by $\mathbb R$ the set of real scalars and by $\mathbb R^{m\times n}$ the set of real $m \times n$ matrices. The transpose and conjugate transpose of a matrix $G$ is denoted by $G^T$ and $G^{\ast}$, respectively. The set $\mathbf S^n$ denotes $n \times n$ Hermitian matrices, and $I_n$ denotes the $n \times n$ identity matrix. We denote by $\mathbf 1_n$ an $n$-dimensional vector with all of its components  equal to~$1$. We denote the real part of a complex vector $v$ by $\real(v)$ and $\eigen(G)$ denotes the eigenvalues of a matrix $G$. We use superscripts for indexing different matrices, and we use subscripts to refer to different elements in the matrix. Hence, by $G^{k}_{ij}$ we denote the element on the $i$th row and $j$th column of the matrix $G^k$. Similarly, $v^k_i$ is the $i$th component of the vector $v^k$. Given matrices $G^k$ for $k = 1, \dots, N$, $\diag(G^1, \dots, G^N)$ denotes a block-diagonal matrix with diagonal blocks specified by the given matrices. Likewise, given vectors $v^k$ for $k= 1, \dots, N$, the column vector $(v^1, \dots, v^N)$ is all of the given vectors stacked. The generalized matrix inequality $G \prec H$ ($G \preceq H$) means that $G-H$ is negative (semi)definite. Given state-space data $A, B, C$ and $D$, we denote by 
\small
\begin{align*}
G(s) := \begin{bmatrix} \begin{array}{c|c} A  & B \\ \hline  C & D  \end{array}\end{bmatrix}
\end{align*}
\normalsize
the transfer function matrix $G(s) = C\left( sI-A \right)^{-1}B + D$. The infinity norm of a transfer function matrix is defined as $\| G(s) \|_\infty = \supremum_{\omega \in \mathbb R} \bar \sigma(G(j\omega))$, where $\supremum$ denotes the supremum and $\bar \sigma$ denotes the largest singular value of a matrix. By $\mathcal L_2^n$ we denote the set of $n$-dimensional square integrable signals, and $\mathcal{RH}_{\infty}^{m \times n}$ represents the set of real, rational $m \times n$ transfer matrices with no poles in the closed right half plane. Given a graph $Q(V,E)$ with vertices $V = \{v_1, \dots, v_n\}$ and edges $E \subseteq V\times V$, two vertices $v_i, v_j \in V$ are adjacent if $(v_i,v_j) \in E$, and we denote the set of adjacent vertices of $v_i$ by $\adj (v_i) = \{ v_j \in V | (v_i, v_j) \in E \}$. The degree of a vertex (node) in a graph is defined as the number of its adjacent vertices, i.e., $\degree(v_i) = |\adj(v_i)|$. The adjacency matrix of a graph $Q(V,E)$ is defined as a $|V| \times |V|$ matrix $A$ where $A_{ij} = 1$ if $(i,j) \in E$ and $A_{ij} =0$ otherwise.
%%%%%%%%%%%%%%%%%%%%%%%%%%%%%%%%%%%%%%%%%%%%
%%%%%%%%%%%%%%%%%%%%%%%%%%%%%%%%%%%%%%%%%%%%
\section{Robustness Analysis Using IQCs} \label{sec:IQC}

%%%%%%%%%%%%%%%%%%%%%%%%%%%%%%%%%%%%%%%%%%%%
\subsection{Integral Quadratic Constraints}
 IQCs play an important role in robustness analysis of uncertain systems where they are used to characterize the uncertainty in the system. Let  $\Delta: \mathbb R^d \rightarrow \mathbb R^d$ denote a bounded and causal operator. This operator is said to satisfy the IQC defined by $\Pi$, i.e., $\Delta \in \IQC(\Pi)$, if 
\small
\begin{align}\label{eq:IQCT}
\int_{0}^{\infty} \begin{bmatrix} v \\ \Delta(v) \end{bmatrix}^{T} \Pi \begin{bmatrix} v \\ \Delta(v) \end{bmatrix} \, dt \geq 0, \quad \forall v \in \mathcal{L}_2^d \ , 
\end{align} 
\normalsize
where $\Pi$ is a bounded and self-adjoint operator. Additionally, the IQC in~\eqref{eq:IQCT} can be written in the frequency domain as
\small
\begin{align}\label{eq:IQCF}
\int_{-\infty}^{\infty} \begin{bmatrix} \widehat{v}(j\omega) \\ \widehat{\Delta(v)}(j\omega) \end{bmatrix}^{\ast} \Pi(j\omega) \begin{bmatrix} \widehat{v}(j\omega) \\ \widehat{\Delta(v)}(j\omega) \end{bmatrix} \, d\omega \geq 0, 
\end{align} 
\normalsize
where $\hat v$ and $\widehat{\Delta(v)}$ are the Fourier transforms of the signals \cite{ulfiqc,rantzer}. IQCs can also be used to describe operators constructed from other operators. For instance, suppose $\Delta^i \in\IQC(\Pi^i)$, where \small$\Pi^i = \begin{bmatrix} \Pi^i_{11} & \Pi^i_{12} \\ \Pi^i_{21} & \Pi^i_{22} \end{bmatrix}$\normalsize.
Then the diagonal operator $\Delta = \diag(\Delta^1, \dots, \Delta^N)$ satisfies the IQC defined by 
\small
\begin{align}\label{eq:IQCDiag}
\bar{\Pi} = \begin{bmatrix} \bar{\Pi}_{11} & \bar{\Pi}_{12}\\ \bar{\Pi}_{21} & \bar{\Pi}_{22}\end{bmatrix},
\end{align}
\normalsize
where $\bar \Pi_{ij} = \diag( \Pi_{ij}^1, \dots,  \Pi_{ij}^N)$, \cite{ulfiqc}. We will use diagonal operators to describe the uncertainties of interconnected systems, but first we briefly describe the IQC analysis framework for robustness analysis of uncertain systems. 
%%%%%%%%%%%%%%%%%%%%%%%%%%%%%%%%%%%%%%%%%%%%
\subsection{IQC Analysis}
Consider the following uncertain system,
\begin{equation}\label{eq:UncertainSystem}
\begin{split}
p = G q, \quad q = \Delta(p),
\end{split}
\end{equation}
where $G \in \mathcal{RH}_{\infty}^{m\times m}$ is referred to as the system transfer function matrix, and $\Delta: \mathbb R^m \rightarrow \mathbb R^m$ is a bounded and causal operator representing the uncertainty in the system. The uncertain system in \eqref{eq:UncertainSystem} is said to be robustly stable if the interconnection between $G$ and $\Delta$ remains stable for all uncertainty values described by $\Delta$, \cite{essentials}. Using IQCs, the following theorem provides a framework for analyzing robust stability of uncertain systems.
\begin{theorem}[IQC analysis, \cite{ulfiqc}]\label{thm:IQC}
The uncertain system in~\eqref{eq:UncertainSystem} is robustly stable, if 
\begin{enumerate}
\item for all $\tau \in [0,1]$ the interconnection described in \eqref{eq:UncertainSystem}, with $\tau\Delta$, is well-posed;
\item for all $\tau \in [0,1]$, $\tau \Delta \in \IQC(\Pi)$;
\item there exists $\epsilon > 0$ such that
\small
\begin{align}\label{eq:thmIQC}
\begin{bmatrix} G(j\omega) \\ I \end{bmatrix}^{\ast} \Pi(j\omega) \begin{bmatrix} G(j\omega) \\ I \end{bmatrix} \preceq -\epsilon I,  \hspace{2mm} \forall \omega \in [0, \infty].
\end{align}
\normalsize
\end{enumerate}
\end{theorem}
\begin{IEEEproof}
See \cite{ulfiqc,rantzer}.
\end{IEEEproof}
The second condition of Theorem~\ref{thm:IQC} mainly imposes structural constraints on $\Pi$, \cite{rantzer}. 
IQC analysis then involves a search for an operator $\Pi$ that satisfies the LMI in \eqref{eq:thmIQC} with the required structure. This can be done using either the KYP lemma-based formulation, \cite{ulfiqc,rantzer,ran:96},  or it can be done approximately by performing frequency gridding where the feasibility of the LMI in~\eqref{eq:thmIQC} is checked only for a finite number of frequencies. In the next section, we propose an efficient method for robust stability analysis of sparsely interconnected uncertain systems based on the frequency-gridding approach.
%%%%%%%%%%%%%%%%%%%%%%%%%%%%%%%%%%%%%%%%%%%%
%%%%%%%%%%%%%%%%%%%%%%%%%%%%%%%%%%%%%%%%%%%%
\section{Robust Stability Analysis of Interconnected Uncertain Systems}\label{sec:RobustStability}

%%%%%%%%%%%%%%%%%%%%%%%%%%%%%%%%%%%%%%%%%%%%
%%%%%%%%%%%%%%%%%%%%%%%%%%%%%%%%%%%%%%%%%%%%
\subsection{Interconnected Uncertain Systems}\label{sec:Inter}
Consider a network of $N$ uncertain subsystems where each of the subsystems is described as
\small
\begin{equation}\label{eq:Subsystems}
\begin{split}
&p^i = G_{pq}^i q^i + G_{pw}^iw^i \\
&z^i = G_{zq}^i q^i + G_{zw}^iw^i\\
&q^i = \Delta^i(p^i),
\end{split}
\end{equation}
\normalsize
where $G_{pq}^i \in \mathcal{RH}_{\infty}^{d_i \times d_i}$, $G_{pw}^i \in \mathcal{RH}_{\infty}^{d_i \times m_i}$, $G_{zq}^i \in \mathcal{RH}_{\infty}^{l_i \times d_i}$, $G_{zw}^i \in \mathcal{RH}_{\infty}^{l_i \times m_i}$, and $\Delta^i:\mathbb{R}^{d_i} \to \mathbb{R}^{d_i}$. If we let $p = (p^1, \dots, p^N)$, $q = (q^1, \dots, q^N)$, $w = (w^1, \dots, w^N)$ and $z = (z^1, \dots, z^N)$, the interconnection among the subsystems in \eqref{eq:Subsystems} can be characterized by the interconnection constraint $w = \Gamma z$ where $\Gamma$ is an interconnection matrix that has the following structure
\small
\begin{align}\label{eq:Interconst}
  \underbrace{ \begin{bmatrix}
    w^1\\w^2\\ \vdots\\ w^N
  \end{bmatrix}}_{w} = 
\underbrace{  \begin{bmatrix}
    \Gamma_{11} & \Gamma_{12} & \cdots & \Gamma_{1N} \\
    \Gamma_{21} & \Gamma_{22} & \cdots & \Gamma_{2N} \\
    \vdots & \vdots & \ddots & \vdots \\
    \Gamma_{N1} & \Gamma_{N2} & \cdots & \Gamma_{NN} 
  \end{bmatrix}}_{\Gamma}
  \underbrace{ \begin{bmatrix}
    z^1\\z^2\\ \vdots\\ z^N
  \end{bmatrix}}_{z}.
\end{align}
\normalsize
Each of the blocks $\Gamma_{ij}$ in the interconnection matrix are 0-1 matrices that describe how individual components of the input-output vectors of different subsystems are connected to each other. The entire interconnected uncertain system can be expressed as
\small
\begin{equation}\label{eq:SysInter}
\begin{split}
p& = G_{pq} q + G_{pw}w \\
z& = G_{zq} q + G_{zw}w\\
q& = \Delta(p)\\
w& = \Gamma z,
\end{split}
\end{equation}
\normalsize
 where $G_{\star\bullet} = \diag(G_{\star\bullet}^1, \dots, G_{\star\bullet}^N)$ and $\Delta  = \diag(\Delta^1, \dots, \Delta^N)$. Using this description of interconnected uncertain systems, we consider two formulations for analyzing robust stability, namely a "lumped" and a "sparse" formulation, as we explain next.
%%%%%%%%%%%%%%%%%%%%%%%%%%%%%%%%%%%%%%%%%%%%
\subsection{Lumped Formulation}\label{sec:Lumped}
The classical approach to robust stability analysis of interconnected uncertain systems is to eliminate the interconnection constraint in \eqref{eq:Interconst} in order to describe the entire interconnected system  as a lumped system 
\begin{equation}\label{eq:Lumped}
\begin{split} 
p = \bar G q, \quad q = \Delta(p),
\end{split}
\end{equation}
where $\bar G = G_{pq} + G_{pw}(I - \Gamma G_{zw})^{-1}\Gamma G_{zq}$. We will refer to $\bar G$ as the lumped system transfer function matrix. Note that $I - \Gamma G_{zw}$ must have a bounded inverse for all frequencies in order for the interconnection to be well-posed. 

Using the lumped formulation, one can use the IQC framework from Theorem \ref{thm:IQC} to analyze the robustness of the interconnected uncertain system. Let $\Delta \in \IQC(\bar \Pi)$ and assume that it satisfies the regularity conditions in Theorem~\ref{thm:IQC}, i.e., conditions 1 and 2 in the theorem. Now suppose that $\bar G \in \mathcal{RH}_{\infty}^{\bar d \times \bar d}$, where $\bar d = \sum_{i=1}^{N}d_i$. The interconnected uncertain system is then robustly stable if there exists a matrix $\bar \Pi$ such that 
\small
\begin{align} \label{eq:IQCLumped}
\begin{bmatrix}
  \bar G(j\omega)\\ I
\end{bmatrix}^\ast \bar \Pi(j\omega)
\begin{bmatrix}
  \bar G(j\omega) \\ I
\end{bmatrix} \preceq -\epsilon I, \quad \forall \omega \in [0, \infty],
\end{align}
\normalsize
for some $\epsilon > 0$. Notice that the matrix on the left hand side of LMI \eqref{eq:IQCLumped} is of order $\bar d$, and it is generally dense even if the subsystems are sparsely interconnected. IQC analysis of large-scale interconnected systems based on the lumped formulation can therefore be prohibitively costly even when the interconnection matrix $\Gamma$ is sparse. This is because the matrix $(I - \Gamma G_{zw})^{-1}$ is generaly dense even if $I - \Gamma G_{zw}$ is sparse. This follows from the Cayley-Hamilton theorem. Next we consider a sparse formulation of the analysis problem. 
%%%%%%%%%%%%%%%%%%%%%%%%%%%%%%%%%%%%%%%%%%%%
\subsection{Sparse Formulation}\label{sec:Sparse}
To avoid solving dense LMIs, we express the interconnection constraint using IQCs. First, notice that the equation $w = \Gamma z$ is equivalent to the following quadratic constraint
\small
\begin{align}\label{eq:IQCSimple}
-\|w-\Gamma z \|_X^2 = -\begin{bmatrix} z \\ w \end{bmatrix}^{\ast}
\begin{bmatrix}
  -\Gamma^T \\ I
\end{bmatrix} X
\begin{bmatrix}
  -\Gamma & I
\end{bmatrix}
\begin{bmatrix} z \\ w \end{bmatrix}\geq 0
\end{align}
\normalsize
where $\| \cdot \|_X$ denotes the norm induced by the inner product $\langle \alpha , X \beta \rangle$ for some $X \succ 0$ of order $\bar m = \sum_{i=1}^N m_i$. The inequality \eqref{eq:IQCSimple} can therefore be rewritten as an IQC defined by the multiplier
\small
\begin{align}
\label{eq:IQCInter} 
\hat{\Pi} &=\begin{bmatrix} -\Gamma^T X \Gamma & \Gamma^T X\\ X\Gamma & -X \end{bmatrix}.
\end{align}
\normalsize
This allows us to include the interconnection constraint in the IQC analysis problem explicitly, and this often results in sparse LMIs if the interconnection matrix is sparse.

Consider \eqref{eq:SysInter}, and let $\Delta \in \IQC(\bar \Pi)$ where $\bar \Pi$ is defined in~\eqref{eq:IQCDiag}. Assuming that  the first and second condition in Theorem \ref{thm:IQC} are satisfied, then it can be shown that the interconnected uncertain system is stable if there exist matrices  $\bar \Pi$ and $X\succ 0$ such that
\small
\begin{multline}\label{eq:IQCInterconnected}
\begin{bmatrix} G_{pq} & G_{pw} \\ I & 0 \end{bmatrix}^{\ast}\begin{bmatrix} \bar{\Pi}_{11} & \bar{\Pi}_{12} \\ \bar{\Pi}_{21} & \bar{\Pi}_{22} \end{bmatrix}\begin{bmatrix} G_{pq} & G_{pw} \\ I & 0 \end{bmatrix} -\\
 \begin{bmatrix} -G_{zq}^{\ast}\Gamma^T\\ I
  -G_{zw}^{\ast}\Gamma^T \end{bmatrix}X\begin{bmatrix} -\Gamma
  G_{zq} & I-\Gamma G_{zw} \end{bmatrix} \preceq -\epsilon I.
\end{multline} 
\normalsize
for $\epsilon > 0$ and for all $\omega \in [0, \infty]$, \cite{and+han:12}. The sparse formulation and lumped formulation result in equivalent optimization problems. The following theorem establishes this equivalence.
\begin{theorem}\label{thm:IQCSparse}
The LMI in \eqref{eq:IQCInterconnected} is feasible if and only if the LMI in~\eqref{eq:IQCLumped} is feasible.
\end{theorem}
\begin{IEEEproof} 
See \cite{and+han:12}. 
\end{IEEEproof}
Theorem \ref{thm:IQCSparse} implies that analyzing robust stability of interconnected uncertain systems using the sparse and lumped formulations lead to equivalent conclusions. As a result, using the sparse formulation of the IQC analysis does not change the conservativeness of the analysis. 
However, note that the LMI in \eqref{eq:IQCInterconnected} is of order $\sum_{i=1}^{N}(d_i + m_i)$, and if the matrix variable $X$ is dense, the LMI will also be dense in general. Hence, if the scaling matrix $X$ is chosen to be dense,  the sparse formulation does not present any improvement compared to the lumped formulation. This issue can be addressed using the following corollary to Theorem \ref{thm:IQCSparse}.

\begin{corollary} \label{cor:IQC}
  In \eqref{eq:IQCInterconnected}, it is sufficient to consider a 
  diagonal scaling matrix of the form $X= x I$ with $x > 0$.
\end{corollary}
\begin{IEEEproof} 
See \cite{and+han:12}. 
\end{IEEEproof}
Corollary \ref{cor:IQC} implies that it is possible to choose the scaling matrix $X$ as a diagonal matrix with a single scalar variable, without adding any conservativeness to the analysis approach. Consequently, if the interconnection matrix is sufficiently sparse, then the LMI in \eqref{eq:IQCInterconnected} will generally also be sparse.
\begin{rem}
The proposed formulation for analyzing interconnected uncertain systems can also be used to analyze systems with uncertain interconnections. This can be done by modifying the subsystems uncertainty descriptions and their system transfer matrices, in order to accommodate the uncertainty in the interconnection within the subsystems uncertainty blocks.
\end{rem}
\begin{rem}
As was mentioned in Section \ref{sec:Intro}, we solve the semi-infinite LMIs in \eqref{eq:IQCLumped} and \eqref{eq:IQCInterconnected} by performing frequency gridding where instead of considering all the frequencies in $[0, \infty]$, we study only a finite set frequencies denoted by $S_f$.   
\end{rem}
%

%%%%%%%%%%%%%%%%%%%%%%%%%%%%%%%%%%%%%%%%%%%%
%%%%%%%%%%%%%%%%%%%%%%%%%%%%%%%%%%%%%%%%%%%%
\section{Sparsity in Semidefinite Programs (SDPs)}\label{sec:Chordal}
In this section, we discuss how the sparse formulation of the robustness analysis problem can be solved efficiently. The approach is based on sparse Cholesky factorization techniques \cite{blp:94}, \cite{geo:93}, which play a fundamental role in many sparse matrix algorithms. Given a sparse positive definite matrix $X$, it is often possible to compute a sparse Cholesky factorization of the form $P^TXP = LL^T$ where $P$ is a permutation matrix and $L$ is a sparse lower-triangular Cholesky factor \cite{geo:81,rot:93,duf:89}. The nonzero pattern of $L+L^T$ depends solely on $P$: it includes all the nonzeros in $P^TXP$ and possibly a number of additional nonzeros which are referred to as \emph{fill}. For general sparsity patterns, the problem of computing a minimum fill factorization is known to be NP-complete \cite{Yan:81}, so in practice a fill-reducing permutation is often used. Sparsity patterns for which a zero-fill permutation exists are called \emph{chordal}, and for such sparsity patterns it is possible to efficiently compute a permutation matrix $P$ that leads to zero fill; see \cite{blp:94} and references therein.

Sparse factorization techniques are also useful in interior-point methods for solving semidefinite optimization problems of the form 
\small
\begin{subequations}\label{eq:SDPDual}
\begin{align}
& \minimize_{S,y} \quad \ b^T y \label{eq:SDPDual1} \\
& \subject \quad  \sum_{i = 1}^{m} y_iQ^i + S  = W, \ \ S \succeq 0. \label{eq:SDPDual2}
 \end{align}
\end{subequations}
\normalsize
The variables are $S \in  \mathbf S^n$ and $y \in \mathbb R^{m}$, and the problem data are $b \in \mathbb R^{m}$ and $W, Q^i \in \mathbf S^n$ for $i = 1, \dots, m$. The LMIs in~\eqref{eq:IQCLumped} and~\eqref{eq:IQCInterconnected} are of the form \eqref{eq:SDPDual} if we let $b = 0$, \cite{boyd:04,elg:00}, and it is clear that the slack variable $S$ inherits its sparsity pattern from the problem data because of the equality constraint \eqref{eq:SDPDual2}. This means that $S$ is sparse if the aggregate sparsity pattern of the data matrices is sparse. 

Solving the problem \eqref{eq:SDPDual} using an interior-point method involves evaluating the logarithmic barrier function $\phi(S) = -\log\det S$, its gradient $\nabla \phi(S) = -S^{-1}$, and terms of the form $\tr(Q^i\nabla^2 \phi(S)Q^j) = \trace(S^{-1}Q^iS^{-1}Q_j)$ at each iteration. When $S$ is sparse, this can be done efficiently using sparse factorization techniques. Note however that $S^{-1}$ is generally not sparse even when $S$ is sparse, but it is possible to avoid forming $S^{-1}$ by working in a lower-dimensional subspace defined by the filled sparsity pattern of $P^TSP$, \cite{and+dah+van10}, \cite{and:12}. This approach generally works well when filled sparsity pattern is sufficiently sparse. The cost of forming and factorizing the so-called Newton equations typically dominates cost of a single interior-point iteration. The Newton equations are of the form
\begin{align}\label{eq:Newton}
H \Delta y = r
\end{align}
where $\Delta y$ is the search direction, $H_{ij} = \trace(S^{-1}Q^i S^{-1} Q^j)$ for $i,j = 1, \dots, m$, and $r$ is a vector of residuals. In general, $H$ is dense even if the data matrices $W$ and $Q^i$ for $i = 1, \dots, m$ are sparse. As a result, it is typically not possible to reduce the computational cost of factorizing $H$, but when the data matrices are sparse and/or have low rank, it is possible to form $H$ very efficiently by exploiting the structure in the data matrices, \cite{ben:99,Fuji:97}.
%%%%%%%%%%%%%%%%%%%%%%%%%%%%%%%%%%%%%%%%%%%%
\section{Numerical Experiments}\label{sec:Results}
In this section, we compare the computational cost of robustness analysis based on the sparse and the lumped formulations using two sets of numerical experiments. In Section~\ref{sec:Chain}, we study a chain of uncertain systems where we compare the performance of the sparse and lumped formulations with respect to the number of subsystems in the network. In Section~\ref{sec:Tree}, we illustrate the effectiveness of the sparse formulation by analyzing an interconnection of uncertain systems over a so-called scale-free network.
%%%%%%%%%%%%%%%%%%%%%%%%%%%%%%%%%%%%%%%%%%%%
\subsection{Chain of Uncertain Systems}\label{sec:Chain}
\begin{figure}
\begin{center}
\begin{tikzpicture}[scale=0.4, inner sep=0pt]

  \tikzstyle{nodeA} = [draw, fill=none, minimum size=1.6cm]
  \tikzstyle{nodeB} = [draw, fill=none, minimum size=1.0cm]

  \draw (0,0) node [style=nodeA]
  {$G^1(s)$};
  \draw (7,0) node [style=nodeA]
  {$G^2(s)$};
  \draw (16.5,0) node [style=nodeA]
  {$G^N(s)$};

  \draw (0,4)  node [style=nodeB]
  {$\delta^1$};
  \draw (7,4)  node [style=nodeB]
  {$\delta^2$};
  \draw (16.5,4) node [style=nodeB]
  {$\delta^N$};
 
  \draw[->] (-2,1) -- (-3,1) -- node[right]{$\phantom{.}p^1$} (-3,4) -- (-1.25,4);
  \draw[->] (1.25,4) -- (3,4) -- node[left]{$q^1\phantom{.}$} (3,1) -- (2,1);
  \draw[->] (2,0) -- node[above]{\raisebox{1.5mm}{$z^1$}} (5,0);
  \draw[->] (5,-1) -- node[below]{\raisebox{-2.5mm}{$z_1^2$}} (2,-1);

  \draw[->] (5,1) -- (4,1) -- node[right]{$\phantom{.}p^2$} (4,4) -- (5.75,4);
  \draw[->] (8.25,4) -- (10,4) -- node[left]{$q^2\phantom{.}$} (10,1) -- (9,1);
  \draw[->] (9,0) -- node[above]{\raisebox{1.5mm}{$z_2^2$}} (11,0);
  \draw[->] (11,-1) --  node[below]{\raisebox{-2.5mm}{$z_1^3$}} (9,-1);

  \draw[->] (14.5,1) -- (13.5,1) -- node[right]{$\phantom{.}p^N$} (13.5,4) -- (15.25,4);
  \draw[->] (17.75,4) -- (19.5,4) -- node[left]{$q^N$} (19.5,1) -- (18.5,1);
  \draw[->] (11.75,0) -- node[above]{\raisebox{1.5mm}{$z_2^{N-1}\hspace{5mm}$}} (14.5,0);
  \draw[->] (14.5,-1) -- node[below]{\raisebox{-2.5mm}{$z^N\hspace{4mm}$}} (11.75,-1);

  \draw (11.5,-0.5) node [fill=none] {{\small  $\cdots$}};

\end{tikzpicture} 
\vspace*{-10pt}
\caption{A chain of $N$ uncertain subsystem.}
\vspace*{-15pt}
\label{fig:System}
\end{center}
\end{figure}
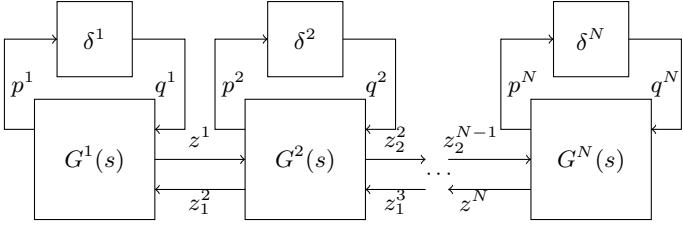

Consider a chain of $N$ uncertain subsystems where each of the subsystems is defined as in \eqref{eq:Subsystems}. We represent the uncertainty in each of the subsystems using scalar uncertain parameters $\delta^1, \dots, \delta ^N$ which correspond to parametric uncertainties in different subsystems. The chain of uncertain systems is shown in Figure \ref{fig:System}. The gains are assumed to be within the normalized interval $[-1, 1]$. The inputs and outputs of the subsystems are denoted by $w^i$ and $z^i$, respectively, where $w^i, z^i \in \mathbb R^2$ for $1 < i < N$, and $w^i, z^i \in \mathbb R$ for $i = 1, N$. 

The interconnections in the network are defined as $w^i_2 = z^{i+1}_1 $ and $w^i_1 = z^{i-1}_2$ for $1<i<N$, and as $w^1 = z^2_1$ and $w^N = z^{N-1}_2$ for the remaining subsystems in the chain, see Figure  \ref{fig:System}. Consequently, the interconnection matrix $\Gamma$ for this network is given by the nonzero blocks \small $\Gamma_{i,i-1} = \Gamma_{i-1,i}^T$ for $i = 2, \dots, N$, where $\Gamma_{i,i-1} = \Gamma_{i-1,i}^T =
  \begin{bmatrix}
    0 & 1 \\ 0 & 0
  \end{bmatrix}, \ i=3,\ldots,N-1$\normalsize, and $ \Gamma_{21} = \Gamma_{12}^T = (1,0), \quad \Gamma_{N-1,N} =  \Gamma_{N,N-1}^T = (0,1)$. Given the uncertainty in each of the subsystems, we have $\delta^i \in \IQC(\Pi^i)$ for $i = 1, \dots, N$, where \small $\Pi^i = \begin{bmatrix} r_i(j\omega) & 0 \\ 0 & -r_i(j\omega) \end{bmatrix}$\normalsize, and $r_i(j\omega) \geq 0$, \cite{rantzer}. We choose the scaling matrix in \eqref{eq:IQCInterconnected} to be of form $X = xI$. Note that analyzing the lumped system yields the LMI in~\eqref{eq:IQCLumped} of order $N$ whereas the sparse LMI in \eqref{eq:IQCInterconnected}
is of order $3N-2$. As a result, for medium-sized networks, it may be computationally cheaper to solve the dense LMI in \eqref{eq:IQCLumped}, but for large and sparse networks, the sparse formulation is generally much more tractable. In order to confirm this, we conduct a set of numerical experiments where we compare the computation time required to solve the lumped and sparse formulation of the analysis problem for different number of subsystems in the chain. The interconnected systems considered in these experiments are chosen such that their robust stability can be established using both sparse and lumped formulations of the analysis problem. Note that generating such systems are generally not straightforward. In this paper, we use the following approach to generate such systems. Consider the interconnected system description in \eqref{eq:Subsystems} and \eqref{eq:SysInter}, and assume that
 \small 
 \begin{align*}
 G_{\star\bullet}^i(s) = \begin{bmatrix} \begin{array}{c|c} A_{\star\bullet}^i  & B_{\star\bullet}^i \\ \hline  C_{\star\bullet}^i & D_{\star\bullet}^i  \end{array}\end{bmatrix},
 \end{align*} 
 \normalsize
 for $i = 1, \dots,N$, $\star \in \{ p,z \}$ and $\bullet \in \{ q,w \}$. This results in 
\small
\begin{align*}
G_{\star\bullet}(s) = \begin{bmatrix} \begin{array}{c|c} A_{\star\bullet}  & B_{\star\bullet} \\ \hline  C_{\star\bullet} & D_{\star\bullet}  \end{array}\end{bmatrix},
\end{align*}
\normalsize
where $A_{\star\bullet} = \diag(A_{\star\bullet}^1, \dots, A_{\star\bullet}^N)$, $B_{\star\bullet} = \diag(B_{\star\bullet}^1, \dots, B_{\star\bullet}^N)$, $C_{\star\bullet} = \diag(C_{\star\bullet}^1, \dots, C_{\star\bullet}^N)$ and $D_{\star\bullet} = \diag(D_{\star\bullet}^, \dots, D_{\star\bullet}^N)$. 
The system transfer matrices for the subsystems are then chosen such that they satisfy the following conditions
\begin{enumerate}
\item $\real(\eigen(A_{\star\bullet}^i)) \prec 0$ for all $i = 1, \dots, N$, $\star \in \{ p,z \}$ and $\bullet \in \{ q,w \}$;
\item there exists $r_i(j\omega) \geq 0$ such that
\small
\begin{align*}
\begin{bmatrix}
   G^i_{pq}(j\omega)\\ I
\end{bmatrix}^\ast  \Pi^i(j\omega)
\begin{bmatrix}
   G^i_{pq}(j\omega) \\ I
\end{bmatrix} \preceq -\epsilon I, \quad \forall \omega \in S_f,
\end{align*}
\normalsize
for $\epsilon > 0$ and all $i = 1, \dots, N$;
\item $\real(\eigen(A_{l})) \prec 0$, where $A_{l} = A_{zw} - B_{zw}(I-\Gamma D_{zw})^{-1}\Gamma C_{zw}$,
\end{enumerate}
where the first and second conditions ensure that each of the subsystems is robustly stable, and the third condition guarantees that the lumped system transfer function matrix satisfies $\bar G \in \mathcal{RH}_{\infty}^{\bar d \times \bar d}$. Also note that the non-singularity of the matrix $I - \Gamma D_{zw}$ is guaranteed by the well-posedness of the interconnection. We generate subsystems that satisfy the three mentioned conditions according to the procedure described in Appendix~\ref{app:sub}. 
We then solve the lumped formulation of the analysis problem to check whether robust stability of the interconnected system can be established using IQC-based analysis methods. If the system is proven to be robustly stable using the lumped formulation, we conduct the analysis once again using its sparse formulation. 

\begin{figure}
\begin{center}
\begin{tikzpicture}
    \begin{semilogyaxis}[font=\scriptsize, xlabel={\small $N$}, ylabel={\small CPU time (seconds)},width=0.85\linewidth,grid=major,ymin=1e-2,ymax=1e3,xmin=0,xmax=200,legend pos=north west,legend style={cells={anchor=west}}]

 \addplot [color=green, mark=*,dashed] coordinates {
(10,2.80e-02)
(20,6.60e-02)
(30,1.57e-01)
(40,3.29e-01)
(50,6.76e-01)
(70,2.18e+00)
(90,4.64e+00)
(100,8.19e+00)
(120,1.46e+01)
(140,2.58e+01)
(160,4.19e+01)
(180,6.32e+01)
(200,1.11e+02) 
}; %SeDuMi

 \addplot [color=red, mark=o] coordinates {
(10,1.33e-01)
(20,2.00e-01)
(30,2.66e-01)
(40,4.00e-01)
(50,7.08e-01)
(70,1.90e+00)
(90,4.43e+00)
(100,6.96e+00)
(120,1.27e+01)
(140,2.16e+01)
(160,3.84e+01)
(180,5.54e+01)
(200,1.01e+02) 
}; %SDPT3

 \addplot [color=blue, mark=square*,dashed] coordinates {
(10,4.80e-02)
(20,1.32e-01)
(30,3.91e-01)
(40,8.21e-01)
(50,1.52e+00)
(70,4.05e+00)
(90,8.79e+00)
(100,1.34e+01)
(120,2.41e+01)
(140,3.90e+01)
(160,5.55e+01)
(180,8.28e+01)
(200,1.18e+02) 
}; %DSDP

 \addplot [color=black, mark=square] coordinates {
 (10,0.1151)
(20,0.2644)
(30, 0.4700)
(40, 0.6982)
(50, 0.9834)
(70, 1.6347)
(90, 2.6257)
(100, 3.2337)
(120, 4.2360)
(140, 5.8132)
(160, 7.0628)
(180, 8.7538)
(200, 12.0270)
 }; %SMCP

\legend{SeDuMi (lumped), SDPT3 (lumped), DSDP (sparse),SMCP (sparse)}
    \end{semilogyaxis}
\end{tikzpicture} 
\caption{Average CPU time required to solve (i) the lumped formulation of the analysis problem
  with SDPT3 and SeDuMi, and (ii) the sparse formulation of the analysis problem with DSDP and SMCP for a chain of length $N$.}
  \vspace*{-18pt}
\label{fig:Results}
\end{center}
\end{figure}
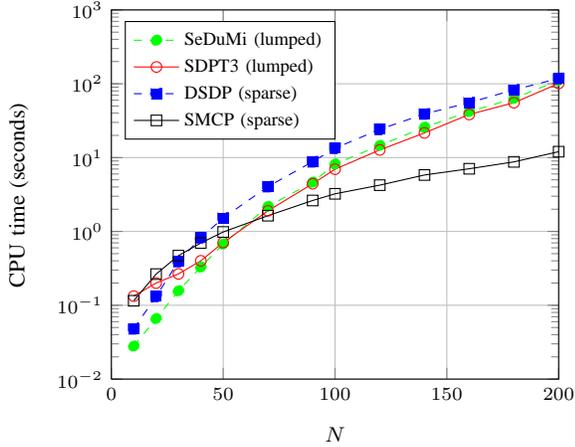

We solve the lumped formulation of the analysis problem using the general-purpose SDP solvers SDPT3, \cite{TTT:99}, and SeDuMi, \cite{Stu:01}. To solve the sparse formulation of the analysis problem, we use the sparse solvers DSDP \cite{BeY:05} and SMCP \cite{and+dah+van10} which exploit the aggregate sparsity pattern. Figure \ref{fig:Results} illustrates the average CPU time required to solve the lumped and sparse formulations of the problem  based on 10 trials for each value of $N$ and for a single frequency. As can be seen from the figure, the analysis based on the sparse LMI in \eqref{eq:IQCInterconnected} is more efficient when the number of subsystems is large, and SMCP yields the best results for $N > 70$. For $N = 200$ the Cholesky factorization of the slack variable $S$ for the sparse formulation was computed separately using the toolbox CHOLMOD, \cite{Chen:08}, with AMD ordering, \cite{Ame:04,Ame:96}, which on average resulted in approximately $1\%$ fill-in. Note that DSDP was not able to match the performance of SMCP in this example. This is because the performance of DSDP mainly relies on exploiting low rank structure in the optimization problem and this structure is not present in this example. We also tried to solve the sparse problem using SDPT3 and SeDuMi, but it resulted in much higher computational times than for the dense problem, so we have chosen not to include these results. 
%
%%%%%%%%%%%%%%%%%%%%%%%%%%%%%%%%%%%%%%%%%%%% 
\vspace*{-5pt}
\subsection{Interconnection of Uncertain Systems Over Scale-free network}\label{sec:Tree}

In this section, we study the interconnection of $N$ uncertain subsystems over a randomly generated scale-free network. Scale-free networks are networks where the degree distribution of the nodes follow a distribution defined as
\small
\begin{align}\label{eq:Plaw}
P_{\alpha}(k) = \frac{k^{-\alpha}}{\sum_{n=1}^{N}n^{-\alpha}}, 
\end{align}
\normalsize
for $k\geq 1$ and where typically $2<\alpha<3$, \cite{clauset:09}. In order to compare the cost of solving the sparse and lumped formulations of the analysis problem, we conduct a set of numerical experiments where we analyze the stability of a network of 500 uncertain subsystems. The network considered in this experiment has been generated using the NetworkX software package, \cite{NetworkX}, which provides us with the adjacency matrix of the network. In this experiment we use $\alpha = 2.5$ and the resulting network is a tree. The degree of each node~$i$ in the network is given by the number of nonzero elements in the $i$th row or column of the adjacency matrix. For this network, the number of nodes with $\degree(i) \leq 5$, $6\leq\degree(i) \leq 10$ and $11 \leq \degree(i)$ are 478, 16 and 6 respectively. Having this degree distribution for the network, we can see that only a small number of nodes in the network have a degree larger than 10. This implies that the interconnection among the subsystems is very sparse. Given the adjacency matrix of the network, we will now describe the interconnection matrix for this network. Recall the interconnection matrix description in \eqref{eq:Interconst}. The only nonzero blocks in the $i$th block-row of this matrix are $\Gamma_{ij}$ where $j \in \adj(i)$. A nonzero block $\Gamma_{ij}$ describes which of the outputs of the $j$th subsystem are connected to what inputs of the $i$th subsystem, and it depends on the chosen indexing for the input-output vectors. More specifically, we use Algorithm~\ref{alg:alg2} to generate the interconnection matrix.
\begin{algorithm}
\caption{Interconnection matrix}\label{alg:alg2}
\footnotesize
\begin{algorithmic}[1] 
\State{Given the adjacency matrix $A$ and the number of subsystems $N$}
\State{Set $\mathbf I = \mathbf 1_N$ and $\mathbf O = \mathbf 1_N$}
\For {$i  = 1, \dots, N$}
\For {$j = 1, \dots, N$}
\If{ $A_{ij}  == 1 $}
\State{Choose $\Gamma_{ij}$ in \eqref{eq:Interconst} such that $w^i_{\mathbf I_i} = z^j_{\mathbf O_j}$}
\State{$\mathbf I_i := \mathbf I_i + 1$, $\mathbf O_j := \mathbf O_j + 1$}
\EndIf
 \EndFor
 \EndFor
 \end{algorithmic}
 \normalsize
 \end{algorithm}
The input-output dimension for the subsystem associated with the $i$th node is given by  $\degree(i)$. In order to provide suitable systems for the experiment in this section, we use the same procedure as described in Section \ref{sec:Chain}, and the uncertainty in the interconnected system is chosen to have the same structure as in Section~\ref{sec:Chain}. 

Table \ref{tab:results} reports the average CPU time required to solve the sparse and lumped formulations of the analysis problem. The results are based on 10 different system transfer matrices for a single frequency, and they clearly demonstrate the advantage of the sparse formulation compared to the lumped formulation. In this experiment, the Cholesky factorization of the slack variable $S$ on average results in about $0.9 \%$ fill-in. Note that the advantage of the sparse formulation becomes even more visible if the number of frequency points considered in the stability analysis, i.e., $|S_f|$, is large. 
\begin{table}
\centering
\caption{\footnotesize time for analyzing 500 subsystems over a scale-free network.\normalsize}
\label{tab:results}
\begin{tabular}{|c||c||c|}
   \hline \textbf{Solver} & \textbf{Avg. CPU time [sec]} & \textbf{Std. dev. [sec]} \\
   \hline
   \hline
   SDPT3 (lumped) &$5640 $ & 529.8 \\
   \hline
   SeDuMi (lumped) & $2760$ &284.3\\
\hline
   DSDP (sparse) & $167 $ & 28.3 \\
   \hline
      SMCP (sparse) & $33 $ & 5.6 \\
      \hline
\end{tabular} 
\end{table}
\section{Conclusions}\label{sec:Conclusions}

IQC-based analysis of sparsely interconnected uncertain systems generally involves solving a set of dense LMIs. In this paper, we have shown that this can be avoided by using IQCs to describe the interconnections among the subsystems. This yields an equivalent formulation of the analysis problem that involves a set of sparse LMIs. By exploiting the sparsity in these LMIs, we have shown that it is often possible to solve the analysis problem associated with large-scale sparsely interconnected uncertain systems more efficiently than with existing methods. As future research directions, we mention the possibility to decompose the sparse LMIs into a set of smaller but coupled LMIs. This would allow us to employ distributed algorithms to solve these LMIs, and such methods may also be amenable to warm-starting techniques to accelerate the frequency sweep in the analysis problem.

\appendices
\section{Subsystems Generation for Numerical Experiments}\label{app:sub}
The system transfer matrices for the subsystems can be generated in different ways. For instance this can either be done by using the \texttt{rss} command in MATLAB$^\text{TM}$ that generates random multi-input multi-output (MIMO) systems, or more generally, by constructing the MIMO systems by randomly choosing each of the elements in the system transfer function matrix separately. The latter approach allows us to produce transfer matrices where different elements have different poles, zeros and orders. In the experiments described in sections~\ref{sec:Chain} and \ref{sec:Tree}, we use the latter approach where we randomly choose each element in the transfer matrices to be a first-order system. Then we explicitly check whether the generated transfer matrices satisfy conditions 1 and 2 in Section \ref{sec:Chain}. 
\begin{algorithm}
\caption{Subsystem rescaling}\label{alg:alg1}
\footnotesize
\begin{algorithmic}[1] 
\State{Given $G^i_{zw}$ for $i= 1, \ldots, N$ and $\gamma$}
\For {$i  = 1, \dots, N$}
\If{ $\| G^i_{zw} \|_\infty \geq 1/\gamma$}
\State{ Choose $\alpha$ such that $\alpha\| G^i_{zw} \|_\infty < 1/\gamma$}
\State{ $G_{zw}^i := \alpha G_{zw}^i$}
\Else { Leave the system transfer function matrix unchanged.}
\EndIf
 \EndFor
 \end{algorithmic}
 \normalsize
 \end{algorithm}

Let $\bar \sigma(\Gamma) = \gamma$ and assume that we have generated system transfer matrices for all subsystems such that they satisfy conditions 1 and~2 in Section \ref{sec:Chain}. The third condition can then be satisfied by scaling the system transfer matrices using Algorithm \ref{alg:alg1}. By the small gain theorem, the rescaled subsystems satisfy the third condition in Section~\ref{sec:Chain}, \cite{essentials}. 
%\section{}

\ifCLASSOPTIONcaptionsoff
  \newpage
\fi
%\IEEEtriggeratref{37}
\bibliographystyle{IEEEtran}  
\bibliography{IEEEabrv,IEEETrans}

\begin{thebibliography}{10}
\providecommand{\url}[1]{#1}
\csname url@samestyle\endcsname
\providecommand{\newblock}{\relax}
\providecommand{\bibinfo}[2]{#2}
\providecommand{\BIBentrySTDinterwordspacing}{\spaceskip=0pt\relax}
\providecommand{\BIBentryALTinterwordstretchfactor}{4}
\providecommand{\BIBentryALTinterwordspacing}{\spaceskip=\fontdimen2\font plus
\BIBentryALTinterwordstretchfactor\fontdimen3\font minus
  \fontdimen4\font\relax}
\providecommand{\BIBforeignlanguage}[2]{{%
\expandafter\ifx\csname l@#1\endcsname\relax
\typeout{** WARNING: IEEEtran.bst: No hyphenation pattern has been}%
\typeout{** loaded for the language `#1'. Using the pattern for}%
\typeout{** the default language instead.}%
\else
\language=\csname l@#1\endcsname
\fi
#2}}
\providecommand{\BIBdecl}{\relax}
\BIBdecl

\bibitem{robustandoptimal}
K.~Zhou, J.~C. Doyle, and K.~Glover, \emph{Robust and Optimal Control}.\hskip
  1em plus 0.5em minus 0.4em\relax Prentice Hall, 1997.

\bibitem{multivariable}
S.~Skogestad and I.~Postlethwaite, \emph{Multivariable Feedback Control}.\hskip
  1em plus 0.5em minus 0.4em\relax Wiley, 2007.

\bibitem{ulfiqc}
U.~J{\"{o}}nsson, ``Lecture notes on integral quadratic constraints,'' May
  2001.

\bibitem{rantzer}
A.~Megretski and A.~Rantzer, ``System analysis via integral quadratic
  constraints,'' \emph{{IEEE} Trans. Autom. Control}, vol.~42, no.~6, pp.
  819--830, Jun. 1997.

\bibitem{ran:96}
A.~Rantzer, ``On the {K}alman-{Y}akubovich-{P}opov lemma,'' \emph{Systems and
  Control Letters}, vol.~28, no.~1, pp. 7--10, 1996.

\bibitem{ChapterBook1}
L.~Vandenberghe, V.~R. Balakrishnan, R.~Wallin, A.~Hansson, and T.~Roh,
  ``Interior-point algorithms for semidefinite programming problems derived
  from the {KYP} lemma,'' in \emph{Positive polynomials in control}, D.~Henrion
  and A.~Garulli, Eds.\hskip 1em plus 0.5em minus 0.4em\relax Springer, Feb
  2005, vol. 312, pp. 195--238.

\bibitem{Anders}
A.~Hansson and L.~Vandenberghe, ``Efficient solution of linear matrix
  inequalities for integral quadratic constraints,'' in \emph{Proceedings of
  the 39th {IEEE} Conference on Decision and Control}, vol.~5, 2000, pp.
  5033--5034.

\bibitem{Parrilo}
P.~A. Parrilo, ``Outer approximation algorithms for {KYP}-based {LMI}s,'' in
  \emph{Proceedings of the American Control Conference}, vol.~4, 2001, pp.
  3025--3028.

\bibitem{Kao2}
C.~Kao, A.~Megretski, and U.~J{\"o}sson, ``Specialized fast algorithms for
  {IQC} feasibility and optimization problems,'' \emph{Automatica}, vol.~40,
  no.~2, pp. 239--252, 2004.

\bibitem{and+dah+van10}
M.~S. Andersen, J.~Dahl, and L.~Vandenberghe, ``Implementation of nonsymmetric
  interior-point methods for linear optimization over sparse matrix cones,''
  \emph{Mathematical Programming Computation}, pp. 1--35, 2010.

\bibitem{fuj+kim+koj+oka+yam09}
K.~Fujisawa, S.~Kim, M.~Kojima, Y.~Okamoto, and M.~Yamashita, ``User's manual
  for {SparseCoLo}: Conversion methods for {SPARSE} {COnnic-form Linear
  Optimization},'' Department of Mathematical and Computing Sciences, Tokyo
  Institute of Technology, Tokyo, Tech. Rep., 2009.

\bibitem{WHJ:09}
R.~Wallin, A.~Hansson, and J.~H. Johansson, ``A structure exploiting
  preprocessor for semidefinite programs derived from the
  kalman-yakubovich-popov lemma,'' \emph{{IEEE} Trans. Autom. Control},
  vol.~54, no.~4, pp. 697--704, Apr. 2009.

\bibitem{Wallin}
R.~Wallin, C.~Kao, and A.~Hansson, ``A cutting plane method for solving
  {KYP-SDP}s,'' \emph{Automatica}, vol.~44, no.~2, pp. 418--429, 2008.

\bibitem{Kao1}
C.~Kao and A.~Megretski, ``A new barrier function for {IQC} optimization
  problems,'' in \emph{Proceedings of the American Control Conference}, vol.~5,
  Jun. 2003, pp. 4281--4286.

\bibitem{Langbort}
C.~Langbort, R.~S. Chandra, and R.~D'Andrea, ``Distributed control design for
  systems interconnected over an arbitrary graph,'' \emph{{IEEE} Trans. Autom.
  Control}, vol.~49, no.~9, pp. 1502--1519, Sep. 2004.

\bibitem{P5}
H.~Fang and P.~J. Antsaklis, ``Distributed control with integral quadratic
  constraints,'' \emph{Proceedings of the 17th {IFAC} World Congress}, vol.~17,
  no.~1, 2008.

\bibitem{UlfLetter}
U.~J{\"{o}}nsson, C.~Kao, and H.~Fujioka, ``A {P}opov criterion for networked
  systems,'' \emph{Systems \& Control Letters}, vol.~56, no. 9--10, pp.
  603--610, 2007.

\bibitem{Kao}
C.~Kao, U.~J{\"o}sson, and H.~Fujioka, ``Characterization of robust stability
  of a class of interconnected systems,'' \emph{Automatica}, vol.~45, no.~1,
  pp. 217--224, 2009.

\bibitem{Ulf}
U.~T. J{\"{o}}nsson and C.-Y. Kao, ``A scalable robust stability criterion for
  systems with heterogeneous {LTI} components,'' \emph{{IEEE} Trans. Autom.
  Control}, vol.~55, no.~10, pp. 2219--2234, Oct. 2010.

\bibitem{Vinnicombe}
I.~Lestas and G.~Vinnicombe, ``Scalable decentralized robust stability
  certificates for networks of interconnected heterogeneous dynamical
  systems,'' \emph{{IEEE} Trans. Autom. Control}, vol.~51, no.~10, pp.
  1613--1625, Oct. 2006.

\bibitem{kim:12}
K.~K. Kim and R.~D. Braatz, ``On the robustness of interconnected or networked
  uncertain linear multi-agent systems,'' \emph{20th International Symposium on
  Mathematical Theory of Networks and Systems}, 2012.

\bibitem{and+han:12}
M.~S. Andersen, A.~Hansson, S.~K. Pakazad, and A.~Rantzer, ``Distributed robust
  stability analysis of interconnected uncertain systems,'' in
  \emph{{Proceedings of the 51st {IEEE} Conference on Decision and Control}},
  2012.

\bibitem{BeY:05}
S.~J. Benson and Y.~Ye, ``{DSDP5} user guide --- software for semidefinite
  programming,'' Mathematics and Computer Science Division, Argonne National
  Laboratory, Argonne, IL, Tech. Rep. ANL/MCS-TM-277, Sep. 2005.

\bibitem{fukuda_exploitingsparsity}
M.~Fukuda, M.~Kojima, , K.~Murota, and K.~Nakata, ``Exploiting sparsity in
  semidefinite programming via matrix completion {I}: General framework,''
  \emph{SIAM Journal on Optimization}, vol.~11, pp. 647--674, 2000.

\bibitem{essentials}
K.~Zhou and J.~C. Doyle, \emph{Essentials of Robust Control}.\hskip 1em plus
  0.5em minus 0.4em\relax Prentice Hall, 1998.

\bibitem{blp:94}
J.~R.~S. Blair and B.~W. Peyton, ``An introduction to chordal graphs and clique
  trees,'' in \emph{Graph Theory and Sparse Matrix Computations}, J.~A. George,
  J.~R. Gilbert, and J.~W.-H. Liu, Eds.\hskip 1em plus 0.5em minus 0.4em\relax
  Springer-Verlag, 1994, vol.~56, pp. 1--27.

\bibitem{geo:93}
A.~George, J.~Gilbert, and J.~Liu, \emph{Graph theory and sparse matrix
  computation}.\hskip 1em plus 0.5em minus 0.4em\relax Springer-Verlag, 1993.

\bibitem{geo:81}
A.~George and J.~Liu, \emph{Computer solution of large sparse positive definite
  systems}.\hskip 1em plus 0.5em minus 0.4em\relax Prentice-Hall, 1981.

\bibitem{rot:93}
E.~Rothberg and A.~Gupta, ``An evaluation of left-looking, right-looking and
  multifrontal approaches to sparse cholesky factorization on
  hierarchical-memory machines.'' \emph{International Journal of High Speed
  Computing}, vol.~5, no.~4, pp. 537--593, 1993.

\bibitem{duf:89}
I.~Duff, A.~Erisman, and J.~Reid, \emph{Direct Methods for Sparse
  Matrices}.\hskip 1em plus 0.5em minus 0.4em\relax Oxford University Press,
  USA, 1989.

\bibitem{Yan:81}
M.~Yannakakis, ``Computing the minimum fill-in is {NP}-complete,'' \emph{SIAM
  J. Algebraic Discrete Methods}, vol.~2, no.~1, pp. 77--79, 1981.

\bibitem{boyd:04}
S.~Boyd and L.~Vandenberghe, \emph{Convex Optimization}.\hskip 1em plus 0.5em
  minus 0.4em\relax Cambridge University Press, 2004.

\bibitem{elg:00}
L.~El~Ghaoui and S.~Niculescu, Eds., \emph{Advances in Linear Matrix Inequality
  Methods in Control}.\hskip 1em plus 0.5em minus 0.4em\relax Society for
  Industrial and Applied Mathematics, 2000.

\bibitem{and:12}
M.~S. Andersen, J.~Dahl, and L.~Vandenberghe, ``Logarithmic barriers for sparse
  matrix cones,'' \emph{Optimization Methods and Software}, 2012.

\bibitem{ben:99}
S.~J. Benson, Y.~Ye, and X.~Zhang, ``Solving large-scale sparse semidefinite
  programs for combinatorial optimization,'' \emph{SIAM Journal on
  Optimization}, vol.~10, no.~2, pp. 443--461, 1999.

\bibitem{Fuji:97}
K.~Fujisawa, M.~Kojima, and K.~Nakata, ``Exploiting sparsity in primal-dual
  interior-point methods for semidefinite programming,'' \emph{Mathematical
  Programming}, vol.~79, pp. 235--253, 1997.

\bibitem{TTT:99}
K.~C. Toh, M.~J. Todd, and R.~H. T{\"u}t{\"u}nc{\"u}, ``{SDPT3} -- {A} {M}atlab
  software package for semidefinite programming, version 1.3,''
  \emph{Optimization Methods and Software}, vol.~11, no.~1, pp. 545--581, 1999.

\bibitem{Stu:01}
J.~F. Sturm, \emph{Using {SEDUMI} 1.02, a Matlab Toolbox for Optimization Over
  Symmetric Cones}, 2001.

\bibitem{Chen:08}
Y.~Chen, T.~A. Davis, W.~W. Hager, and S.~Rajamanickam, ``Algorithm 887:
  {CHOLMOD}, supernodal sparse {Cholesky} factorization and update/downdate,''
  \emph{{ACM} Transactions on Mathematical Software}, vol.~35, no.~3, Oct.
  2008.

\bibitem{Ame:04}
P.~R. Amestoy, T.~A. Davis, and I.~S. Duff, ``Algorithm 837: {AMD}, an
  approximate minimum degree ordering algorithm,'' \emph{{ACM} Transactions on
  Mathematical Software}, vol.~30, no.~3, Sep. 2004.

\bibitem{Ame:96}
------, ``An approximate minimum degree ordering algorithm,'' \emph{{SIAM}
  Journal of Matrix Analysis and Applications}, vol.~17, no.~4, Oct. 1996.

\bibitem{clauset:09}
A.~Clauset, C.~R. Shalizi, and M.~E.~J. Newman, ``Power-law distributions in
  empirical data,'' \emph{SIAM Review}, vol.~51, no.~4, pp. 661--703, 2009.

\bibitem{NetworkX}
A.~A. Hagberg, D.~A. Schult, and P.~J. Swart, ``Exploring network structure,
  dynamics, and function using {NetworkX},'' in \emph{Proceedings of the 7th
  Python in Science Conference}, 2008, pp. 11--15.

\end{thebibliography}

\end{document}